\documentclass{article}
\usepackage{tikz}
\usepackage{subfigure}
\usepackage[english]{babel}
\usepackage{amsmath}
\usepackage{amsthm}
\usepackage{amssymb}
\usepackage{mathrsfs}
\usepackage{amscd}
\usepackage{graphicx}

\def\supp{\mbox{supp}}

\def\({\left(}
\def\){\right)}

\newtheorem{thm}{Theorem}[section]
\newtheorem{cor}[thm]{Corollary}
\newtheorem{lem}[thm]{Lemma}

\newtheorem{prob}[thm]{Problem}

\newtheorem{obs}[thm]{Observastion}
\begin{document}
\title{ The maximal sum of sizes of cross intersecting families for multisets}
\author{Hongkui Wang$^a$, \quad Xinmin Hou$^{a,b}$\footnote{Email: xmhou@ustc.edu.cn (X. Hou), whk98@mail.ustc.edu.cn (H. Wang)}\\
\small $^{a}$ School of Mathematical Sciences\\
\small University of Science and Technology of China, Hefei, Anhui 230026, China.\\
\small$^b$ Hefei National Laboratory\\
\small University of Science and Technology of China, Hefei 230088, Anhui, China
}
\date{}
\maketitle

\begin{abstract}
Let $k$, $t$ and $m$ be positive integers. 
A $k$-multiset of $[m]$ is a collection of $k$ elements of $[m]$ with repetition and without ordering. 
We use $\left(\binom {[m]}{k}\right)$ to denote all the $k$-multisets of $[m]$.
Two multiset families $\mathcal{F}$ and $\mathcal{G}$ in $\left(\binom {[m]}{k}\right)$ are called cross $t$-intersecting if $|F\cap G|\geq t$ for any $F\in \mathcal{F}$ and $G\in \mathcal{G}$. Moreover, if $\mathcal{F}=\mathcal{G}$, we call $\mathcal{F}$ a $t$-intersecting family in $\left(\binom {[m]}{k}\right)$.
Meagher and Purdy~(2011) presented a multiset variant of Erd\H{o}s-Ko-Rado Theorem for $t$-intersecting family in $\left(\binom {[m]}{k}\right)$ when $t=1$, and F\"uredi, Gerbner and Vizer~(2016) extended this result to general $t\ge 2$ with $m\geq 2k-t$, verified a conjecture proposed by Meagher and Purdy~(2011).
In this paper, we determine the maximum sum of cross $t$-intersecting families $\mathcal{F}$ and $\mathcal{G}$ in $\left(\binom {[m]}{k}\right)$ and characterize  the extremal families achieving the upper bound.
For $t=1$ and $m\geq k+1$, the method involves constructing a bijection between multiset family and set family while preserving the intersecting relation.  
 For $t\ge 2$ and $m\ge 2k-t$, we employ  a shifting operation,  specifically the down-compression, which was initiated by F\"uredi, Gerbner and Vizer~(2016).
These results extend the sum-type intersecting theorem for set families originally given by Hilton and Milner (1967).
	
\end{abstract}

\section{Introduction}
We write $[n]$ for $\{1,2,...,n\}$ and $2^{[n]}$ for the power set of $[n]$.  An $\mathcal{F} \subset 2^{[n]}$ is called a {\em t-intersecting family} if the intersection of every two sets in $\mathcal{F}$ contains at least $t$ elements, where $t$ is an integer. We say an 1-intersecting family as intersecting family for convenience.  Let $\binom {[n]}{k}$ denote the family of all $k$-subsets of $[n]$. A set family $\mathcal{F}$ is called $k$-uniform if $\mathcal {F}\subset \binom {[n]}{k}$.
The famous Erd\H{o}s-Ko-Rado Theorem  states that:
\begin{thm}[Erd\H{o}s-Ko-Rado Theorem, \cite{EKR}]\label{THM: EKR}
	$|\mathcal {F} |\leq \binom{n-1}{k-1}$ if $\mathcal{F}$ is a $k$-uniform intersecting family for $n\geq 2k$. Moreover,  the equality holds when $n>2k$ if and only if $\mathcal{F}=\{F\in \binom {[n]}{k} : i\in F\}$ for some $i\in[n]$.
\end{thm}
For the general $t$-intersecting case, Erd\H{o}s, Ko, and Rado \cite{EKR} showed that if $\mathcal{F}\subset \binom{[n]}{k}$ and $\mathcal{F}$ is $t$-intersecting with $n\geq n_0(k,t)$, then $|\mathcal{F}|\leq \binom{n-t}{k-t}$. 
 Frankl \cite{FR} and, subsequently  Wilson~\cite{RMW}, proved that this upper bound is also valid for $n\ge (t+1)(k-t+1)$. Frankl and F\"uredi~\cite{AU}, and late Ahlswede and Khachatrian~\cite{AK}, completely determined the extremal values of $t$-intersecting families for $2k-t\leq n\leq (t+1)(k-t+1)$.

There are number of extensions of Theorem~\ref{THM: EKR}, see a survey \cite{PN}. One important variant of intersecting families is the cross intersecting families. Given positive integers $n$, $k$, $l$, $t$, 
Given a positive integer $t$, two families $\mathcal{A}, \mathcal{B} \in \binom {[n]}{k}$ are {\em cross $t$-intersecting} if  $|A \cap B| \geq t$ for any $A\in \mathcal{A}$ and $B\in \mathcal{B}$. Such a pair $(\mathcal{A},\mathcal{B})$ is called a {\em cross $t$-intersecting pair}. In the rest of this paper, we do not distinguish cross intersecting and cross $1$-intersecting.

For cross intersecting families, Hilton and Milner~\cite{HM} showed the following theorem.   
\begin{thm}[Hilton-Milner Theorem~\cite{HM}]\label{THM:HM}
  Let $n$, $k$ be positive integers with $\mathcal{A} ,\mathcal{B}\subset \binom{[n]}{k}$. Let $\mathcal{A}$ and $\mathcal{B}$ be non-empty cross-intersecting families with $n\geq 2k$. Then   $$|\mathcal{A}|+|\mathcal{B}|\leq \binom{n}{k}-\binom{n-k}{k}+1$$
\end{thm}
Frankl and Tokushige~\cite{FT} generalized the result of Hiton and Milner to set families 
for $\mathcal{A}\subset \binom{[n]}{k}$ and $\mathcal{B}\subset \binom{[n]}{l}$.
\begin{thm}[\cite{FT}]\label{THM:FT}
   Let $n$, $k$, $l$ be positive integers with $\mathcal{A}\subset \binom{[n]}{k}$ and $\mathcal{B}\subset \binom{[n]}{l}$. If $\mathcal{A}$ and $\mathcal{B}$ are non-empty cross intersecting families with $n\geq k+l$ and $k\leq l$. Then   $$|\mathcal{A}|+|\mathcal{B}|\leq \binom{n}{l}-\binom{n-k}{l}+1$$
\end{thm}
Wang and Zhang~\cite{JH} extended Theorems~\ref{THM:HM} and \ref{THM:FT} to cross $t$-intersecting families and characterized the extremal families.
\begin{thm}[\cite{JH}]\label{THM:JH}
	Let $n$, $k$, $l$, $t$ be positive integers with $n\geq 4$, $k, l\geq 2$, $t< \min\{k,l\}$, $k+l< n+t$, $(n,t)\neq (k+l,1)$, and $\binom{n}{k}\leq \binom{n}{l}$. If $\mathcal{A}\subset \binom{[n]}{k}$ and $\mathcal{B}\subset \binom{[n]}{l}$ are non-empty cross $t$-intersecting families, then $$|\mathcal{A}|+|\mathcal{B}|\leq \binom{n}{l}-\sum\limits_{i=0}^{t-1}\binom{k}{i} \binom{n-k}{l-i}+1.$$
Moreover,

 (1) when $\binom{n}{k}\textless \binom{n}{l}$, equality holds if and only if $\mathcal{A}=\{A\}$ and  $\mathcal{B}=\binom {[n]}{l}\backslash N(A)$ for any $A\in \binom {[n]}{k}$, where $N(A)=\{B\, :\, B\in \mathcal{B} \text{ and } B\cap A\neq \emptyset\}$.

(2) when $\binom{n}{k}= \binom{n}{l}$, equality holds if and only if either $\mathcal{A}=\{A\}$ and $\mathcal{B}=\binom {[n]}{l}\backslash N(A)$ for any $A\in \binom {[n]}{k}$, or  $\mathcal{B}=\{B\}$ and  $\mathcal{A}=\binom {[n]}{A}\backslash N(B)$ for any $B\in \binom {[n]}{l}$, or $\{k,l,t\}=\{2,2,1\}$ and $\mathcal{A}=\mathcal{B}=\{C\in \binom {[n]}{2}: i\in C\}$ for some $i\in [n]$, or $\{k,l,t\}=\{n-2,n-2,n-3\}$ and $\mathcal{A}=\mathcal{B}=\binom{A}{n-2}$ for some $A\in \binom{[n]}{n-1}$, where $N(B)=\{A\, :\, A\in \mathcal{A} \text{ and } A\cap B\neq \emptyset\}$.
\end{thm}
There are also several other extensions of Theorems~\ref{THM:HM} and \ref{THM:FT}, see \cite{PF}, \cite{PB}, \cite{PC}. 

In this paper, we concern the intersecting families problem in multisets. We think of a $k$-multiset of $[m]$ as choosing $k$ elements of $[m]$ with repetition and without ordering. We use $\left(\binom {[m]}{k}\right)$ to denote all the $k$-multisets of $[m]$.
 Let $m(i,F)$ denote the multiplicity of $i$ within the multiset $F$, which represents the number of times $i$ appears in $F$.
The cardinality of $F$ is defined as the sum of multiplicities of all elements in $F$, i.e.,  $|F|=\sum_{i=1}^{m}m(i,F)$.
Conventionally, multisets are represented by listing their elements along with their respective multiplicities, formulated as $F=\{1* m(1,F), 2*m(2,F),\dots, m* m(m, F)\}$. 
The {\em support} of a multiset $F$, denoted as supp($F$), in $\left(\binom {[m]}{k}\right)$ is the set of elements $i$ in $[m]$ with $m(i, F)>0$, i.e., supp($F$)$=\{ i\in[m]\, : \, m(i,F)>0\}$. For a multiset family $\mathcal{F}\in \left(\binom {[m]}{k}\right)$, define $\supp(\mathcal{F})=\{supp(F)\,:\,F\in \mathcal{F}\}$.
The intersection of two multisets $F$ and $G$ is defined as a multiset where each element occurs with a multiplicity equal to the minimum of its occurrences in both mutisets, i.e., 
$F\cap G=\{1* \min\{m(1,F),m(1,G)\},..., m*\min\{m(m,F),m(m,G)\}\}$.
A multiset family $\mathcal{F} \subset \left(\binom {[m]}{k}\right)$ is called {\em $t$-intersecting} if $|F_1\cap F_2|\geq t$ for any two multisets $F_1,F_2\in \mathcal{F}$.
We refer a 1-intersecting family as an intersecting family. 
Two multiset families $\mathcal{A}, \mathcal{B}\subset \left(\binom {[m]}{k}\right)$ are said to be {\em cross $t$-intersecting} if  for every pair $A\in \mathcal{A}$ and $B\in \mathcal{B}$, $|A \cap B| \ge t$. We refer a cross 1-intersecting family as a cross intersecting family. We also refer to $(\mathcal{A},\mathcal{B})$ as a {\em cross $t$-intersecting pair}. 
For intersecting families in multisets, Meagher and Purdy~\cite{MP} proved the following result, as partial answer to a conjecture proposed by Brockman and Kay~\cite{BK08}. 
\begin{thm}[\cite{MP}]\label{THM:PM}
 If $m\geq k+1$ and $\mathcal{F}\in \left(\binom {[m]}{k}\right)$ is intersecting, then 
$$|\mathcal{F}|\leq \binom{m+k-2}{k-1}.$$
If $m>k+1$, then equality holds if and only if $\mathcal{F}=\left\{F\in\left(\binom {[m]}{k}\right)\,:\,i\in F\right\}$ for some fixed $i\in [m]$.
\end{thm}
They further conjectured that for $m\ge t(k-t)+2$ and $1\le t\le k$, if  $\mathcal{F}$ is a $t$-intersecting family, then $\mathcal{F}\leq \binom{m+k-t-1}{k-t}$.
F\"uredi, Gerbner and Vizer~\cite{FZ} completely verified this conjecture.
 Meagher and Purdy~\cite{KA}  characterized the extremal structures of the nontrivial intersecting family with $m\ge k+1$, and $t$-intersecting families with $m\geq 2k-t$ later.



In this paper, we focus on the maximum sum of cross $t$-intersecting families for multisets. We first give a generalization of Hilton-Milnor Theorem to multisets.  
Let $\mathcal{F},  \mathcal{G}\in \left(\binom{[m]}{k}\right)$ be  non-empty cross intersecting families. By symmetry, we always assume $|\mathcal{F}|\leq |\mathcal{G}|$. 
\begin{thm}\label{THM: 1-intersecting}
	Let $m\geq k+1$ and let $\mathcal{F},  \mathcal{G}\in \left(\binom{[m]}{k}\right)$ be  non-empty cross intersecting families. Then $$|\mathcal{F}|+|\mathcal{G}|\leq 1+\binom{m+k-1}{k}-\binom{m-1}{k}.$$
Moreover, equality holds if and only if either for $k>2$, $\mathcal{F}=\{[k]\}$ and $\mathcal{G}=\{G\in (\binom{[m]}{k})\,:\,G\cap[k]\neq \emptyset\}$ up to isomorphism, or for $k=2$, $m>3$, $\mathcal{F}=\mathcal{G}=\{D\in (\binom {[m]}{2}): i\in D\}$ or $\mathcal{F}=\{[2]\}$ and $\mathcal{G}=\{G\in (\binom{[m]}{k})\,:\,G\cap[2]\neq \emptyset\}$
up to isomorphism,  or for $\{m,k\}=\{3,2\}$, $\mathcal{F}=\mathcal{G}=\{(12),(23),(13)\}$, or $\mathcal{F}=\{(12)\}$ and $\mathcal{G}=\{(11),(12),(13),(22),(23)\}$, or $\mathcal{F}=\mathcal{G}=\{(11),(12),(13)\}$ up to isomorphism.
\end{thm}

For cross $t$-intersecting families with $2\le t\le k$, we have the following result. 
\begin{thm}\label{THM: t-intersecting}
Let $2\leq t\leq k$ and  $m\geq 2k-t$. If $\mathcal{F}$, $\mathcal{G}\in \left(\binom {[m]}{k}\right)$ are non-empty cross $t$-intersecting families, then 
$$|\mathcal{F}|+|\mathcal{G}|\leq \binom{m+k-1}{k}-\sum\limits_{i=0}^{t-1}\binom{k}{i} \binom{m-1}{k-i}+1.$$
Moreover,  equality holds if and only if $\mathcal{F}=\{[k]\}$, and $\mathcal{G}=\{G\in(\binom{m}{k})\,:\, G\cap[k]\geq t\}$ up to isomorphism.
\end{thm}



The rest of this article is arranged as follows. We give some notations and preliminaries in Section 2. The proof of Theorem~\ref{THM: 1-intersecting} will be given in Section 3. The proof of Theorem~\ref{THM: t-intersecting} will be given in Section 4. We conclude with some discussions and remarks.

\section{Notations and Preliminaries}
A multiset family pair $(\mathcal{F},\mathcal{G})$ is sum maximal cross $t$-intersecting if $\mathcal{F}$ and $\mathcal{G}$ are cross $t$-intersecting multiset families, but adding any other multiset to $\mathcal{F}$ or $\mathcal{G}$ would destroy the cross $t$-intersecting property.
 
We use a special representaion of multisets, which we will employ extensively in the following results. Let $m$ and $l$ be positive integers, and let $M(m,l)=\{(i,j) : 1\leq i\leq m,1\leq j\leq l\}$ be a $m\times l$ rectangle with $m$ rows and $l$ columns. We say $A$ is a subset of $M(m,l)$ if for any $(i,j)\in A$, $A$ follows that $(i,j')\in A$ for all $j'\leq j$. We refer to a subset $A\subseteq M(n,l)$ as a $k$-multiset if the cardinality of $A$ is $k$. It is clear that a subset $A\subseteq M(m,l)$ corresponding to a multiset $A\in\left(\binom{[m]}{k}\right)$, where the multiplicity of the element $i$ equals to $m(i,A)=\max\{s : (i,s)\in A\}$. In the following of this article, we do not distinguish a multiset $A\in \left(\binom{[m]}{k}\right)$ and its corresponding set $A\subseteq M(m,l)$. Specially , $[m]$ is $M(m,1)$.

We need the uniform version of Theorems~\ref{THM:JH} in the proofs of our main theorems.
\begin{cor}[\cite{JH}]\label{COR:JH}
Let $n$, $k$, $t$ be positive integers with $k>t\geq 1$, $n>2k-t$, and $(n,t)\neq (2k,1)$. If $\mathcal{A}, \mathcal{B}\subset \binom{[n]}{k}$ are cross $t$-intersecting families, then $$|\mathcal{A}|+|\mathcal{B}|\leq \binom{n}{k}-\sum\limits_{i=0}^{t-1}\binom{k}{i} \binom{n-k}{k-i}+1.$$
Moreover,
 equality holds if and only if either $\mathcal{A}=\{[k]\}$ and $\mathcal{B}=\{B \,:\, B\cap [k]\not=\emptyset\}$ up to isomorphism, or $(k,t)=(2,1)$ and $\mathcal{A}=\mathcal{B}=\{C\in \binom {[n]}{2} : i\in C\}$ for some fixed $i\in [n]$, or $(k,t)=(n-2,n-3)$ and $\mathcal{A}=\mathcal{B}=\binom{A}{n-2}$ for some $A\in \binom{[n]}{n-1}$.
\end{cor}
For  cross-intersecting families, the following is a simple observation. 
\begin{obs}\label{OBS: complment}
Let $(\mathcal{F},\mathcal{G})$ be a cross $t$-intersecting multiset family pair  on the ground set $[n]$.    
For every $F\in\mathcal{F}$,  the subset of  $[n] \setminus {\supp(F)}$ can not be contained in any multiset $G\in\mathcal{G}$.
\end{obs}
 
\section{Proof of Theorem~\ref{THM: 1-intersecting}}
The key method involves constructing a bijection between multiset family and set family which preserving the intersecting relation. 
Let  $f$ be a map from multiset families $\left(\binom{[n]}{k}\right)$ to $\left(\binom{[m]}{k}\right)$ and let $\mathcal{F}\in \binom{[n]}{k}$. Denote $f(\mathcal{F})=\{f(F)\,:\,F\in \mathcal{F}\}$. 
We have the following proposition.


\begin{lem}\label{Prop:MK}
Let $m$ and $k$ be positive integers, and let $n=m+k-1$. Then there exists a map $f : \binom{[n]}{k} \rightarrow  \left(\binom{[m]}{k}\right)$	with the following properties:
  
 (1) $f$ is a bijection.
  
 (2) For any $B\in \binom{[n]}{k}$, the support of $f(B)$ equals $B\cap [m]$.
  
 (3) If $\mathcal{A}, \mathcal{B}\subseteq \left(\binom{[m]}{k}\right)$ are non-empty cross-intersecting families, 
 then $f^{-1}(\mathcal{A})$ and $f^{-1}(\mathcal{B})$ are  cross-intersecting families in $\binom{[n]}{k}$.
\end{lem}
\begin{proof}
Fix a subset $A\subseteq [m]$ with cardinality $a$, where $1\leq a \leq k$. On one hand, let $\mathcal{B}=\{B\,|\, B \in \binom{[n]}{k} \text{ with } B\cap[m]=A\}$. Then $|\mathcal{B}|=\binom{n-m}{k-a}=\binom{k-1}{k-a}=\binom{a+(k-a)-1}{k-a}=\left|\left(\binom{[a]}{k-a}\right)\right|$.
On the other hand, the number of multisets in $\left(\binom{[m]}{k}\right)$ with support set $A$ is $\binom{k-1}{k-a}$ (which equals to the number of the non-negative solutions of the  indeterminate equation $x_1+x_2+\dots+x_a=k-a$, where $x_1,\dots,x_a$ are non-negetive integers).
Therefore, we can establish a bijection $f : \binom{[n]}{k} \rightarrow  \left(\binom{[m]}{k}\right)$ such that for any set $B\in \binom{[n]}{k}$, the image $f(B)$ is a multiset in $\left(\binom{[m]}{k}\right)$ with support $B\cap [m]$, where each vertex $i\in B\cap [m]$ is assigned a multiplicity of $x_i+1$, and $(x_1, x_2, \ldots, x_{|B\cap[m]|})$ is a non-negative solution of the  indeterminate equation $x_1+x_2+\dots+x_{|B\cap [m]|}=k-|B\cap [m]|$.

Now assume $\mathcal{A},  \mathcal{B}\subseteq \left(\binom{[m]}{k}\right)$ are non-empty cross-intersecting families. 
Then $f^{-1}(\mathcal{A}), f^{-1}(\mathcal{B})\in \binom{[n]}{k}$.
Choose $A\in f^{-1}(\mathcal{A})$ and $B\in f^{-1}(\mathcal{B})$. Then $f(A)\in \mathcal{A}$ and $f(B)\in \mathcal{B}$. Thus $f(A)\cap f(B)\neq \emptyset$, i.e. $\supp(f(A)\cap f(B))\neq \emptyset$. According to (2), $A\cap[m]=\supp(f(A))$ and $B\cap[m]=\supp(f(B))$. Hence $|A\cap B|\geq |A\cap B\cap [m]|$=$|(A\cap[m])\cap (B\cap [m])|=|\supp(f(A))\cap \supp(f(B))|$=$|\supp(f(A)\cap f(B))|>0$.
Therefore,  $f^{-1}(\mathcal{A})$ and $f^{-1}(\mathcal{B})$ are  cross-intersecting families in $\binom{[n]}{k}$.
\end{proof}

Now we are ready to prove Theorem~\ref{THM: 1-intersecting}.
\begin{proof}[Proof of Theorem~\ref{THM: 1-intersecting}:]
Set $n=m+k-1$. According to Lemma~\ref{Prop:MK}, there exists a bijection $f : \binom{[n]}{k} \rightarrow  \left(\binom{[m]}{k}\right)$ that satisfies the properties (2) and (3) in Lemma~\ref{Prop:MK}. 
Let $(\mathcal{F},\mathcal{G})$ be a  sum maximal cross-intersecting pair in $\left(\binom{[m]}{k}\right)$. 
By Theorem~\ref{THM:HM}, we have 
\begin{eqnarray*}
 |\mathcal{F}|+|\mathcal{G}|&=&|f^{-1}(\mathcal{F})|+|f^{-1}(\mathcal{G})|\\
                            &\le&1+\binom{n}{k}-\binom{n-k}{k}\\
                            &=& 1+\binom{m+k-1}{k}-\binom{m-1}{k}.   
\end{eqnarray*}

If $k>2$, by Corollary~\ref{COR:JH}, 
 equality holds if and only if $f^{-1}(\mathcal{F})=\{[k]\}$ and $f^{-1}(\mathcal{G})=\{G'\in \binom{[m+k-1]}{k} \,:\, G'\cap [k] \neq \emptyset\}$ up to isomorphism.
Since $f$ is a bijection from $\binom{[n]}{k}$ to $\left(\binom{[m]}{k}\right)$, we have $\supp(f([k]))=[k]\cap[m]=[k]$ according to (2) of Lemma~\ref{Prop:MK}. Hence $\mathcal{F}=\{[k]\}$. For any $G'\in f^{-1}(\mathcal{G})$, according to (2) of Proposition~\ref{Prop:MK}, $\supp(f(G'))=G'\cap[m]\supseteq G'\cap[k]\not=\emptyset$. Hence $f(G')\cap[k]\not=\emptyset$. Therefore, $\mathcal{G}=\{f(G') \, :\, G'\in f^{-1}(\mathcal{G})\}=\{ G \, :\, G\cap[k]\not=\emptyset\}$.

When $(k,1)=(2,1)$ and $n=m+k-1\not=2k$, i.e., $k=2$ and $m>3$. According to Corollary~\ref{COR:JH}, except for the extremal structure $\mathcal{F}=\{[2]\}$ and $\mathcal{G}=\{G\in \left(\binom{[m]}{2}\right)\,:\,G \cap[2]\neq \emptyset\}$, we also have $f^{-1}(\mathcal{F})$= $f^{-1}(\mathcal{G})$=$\{C\in \binom {[m+1]}{2}: 1\in C\}$ up to isomorphism. For the latter case, we have $1\in\supp(f(C))$ for each $C\in f^{-1}(\mathcal{F})=f^{-1}(\mathcal{G})$. Thus $\supp(\mathcal{F})=\supp(\mathcal{G})=\{\supp(f(C))\, :\, 1\in\supp(f(C))\}\subseteq\{(1),(12),\dots,(1m)\}$, where we denote $(ij)$ as the set $\{i,j\}$ in $\left(\binom{[m]}{2}\right)$.  Note that $|\mathcal{F}|=|\mathcal{G}|=|f^{-1}(\mathcal{G})|$=$m$. Therefore, $\mathcal{F}=\mathcal{G}= \{(11),(12),\dots (1m)\}$.

When $(k,1)=(n-2,n-3)=(m+1-2, m+1-3)$, we have $1=m+1-3$ and $k=m+1-2=2$, i.e., $m=3$ and $m+1=4=2k$.  Corollary~\ref{COR:JH} does not applicable for this case.  At this situation, we have $|\mathcal{F}|+|\mathcal{G}|$ = $1+\binom{4}{2}-\binom{2}{2}=6$. Thus we can list all possible cases for this case. 

If $|\mathcal{F}|$=$1$, then $\mathcal{F}=\{(11)\}$ or $\{(12)\}$ up to isomorphism. If $\mathcal{F}=\{(11)\}$, then $\mathcal{G}\subseteq\{(11),(12),(13)\}$. This is impossible since $|\mathcal{F}|+|\mathcal{G}|=6$.  Hence $\mathcal{F}=\{(12)\}$ and $\mathcal{G}\subseteq\{(11),(12),(13),(22),(23)\}$.
Since $|\mathcal{F}|+|\mathcal{G}|=6$, we have $\mathcal{G}=\{(11),(12),(13),(22),(23)\}$.

If $|\mathcal{F}|$=$2$, then $\mathcal{F}=\{(12),(33)\}$ or $\mathcal{F}=\{(12),(13)\}$, up to isomorphism. For the former case, we have $\mathcal{G}\subseteq\{(13), (23)\}$. This contradicts to $|\mathcal{F}|+|\mathcal{G}|=6$ again. For the latter case, $\mathcal{G}\subseteq\{(11), (12), (13), (23)\}$. Therefore, $\mathcal{G}=\{(11), (12), (13), (23)\}$ since $|\mathcal{F}|+|\mathcal{G}|=6$.

If $|\mathcal{F}|$=$3$, we claim that  $\mathcal{F}$ must be an intersecting family in this case. Suppose not, there exist $F_1$ and $F_2$ in $\mathcal{F}$ such that $F_1\cap F_2$=$\emptyset$. Note that $\left(\binom{[3]}{2}\right)=\{(12), (13), (23), (11), (22), (33)\}$. We define $(12)^{c}=(33)$, $(13)^{c}=(22)$, $(23)^{c}=(11)$, $(33)^{c}=(12)$, $(22)^{c}={13}$ and $(11)^{c}=(23)$. Clearly, $F_1, F_2\notin \mathcal{G}$, otherwise, it will destroy the property of cross-intersecting. According to Observation~\ref{OBS: complment}, we have $F^{c}_1, F^{c}_2\notin \mathcal{G}$.
Thus $\mathcal{G}$ has at most 2 elements, a contradiction to $|\mathcal{F}|+|\mathcal{G}|=6$. Therefore, we have that both $\mathcal{F}$ and $\mathcal{G}$ are intersecting family. Since $(\mathcal{F}, \mathcal{G})$ is a cross intersecting pair, we must have that either $\mathcal{F}=\mathcal{G}=\{(11), (12), (13)\}$ or $\mathcal{F}=\mathcal{G}=\{(12), (13), (23)\}$.

\end{proof}

\section{Proof of Theorem~\ref{THM: t-intersecting}}

 
 
Shifting operations are frequently employed to tackle intersecting family problems. In this section, we use a special shifting operation which is inspired by F\"uredi, Gerbner and Vizer~\cite{FZ} (down-compression).

For a multiset $F\in \left(\binom{[m]}{k}\right)$ with $|F|=k$ and two integers $i,j\in [m]$ and an integer $s\leq m(i,F)$, define a multiset $F'$ with respect to $(F, i,j,s)$ as $\left(F\setminus\cup_{t=s}^{m(i,F)}\{(i,t)\}\right) \bigcup \left(\cup_{l=1}^{m(i,F)-s+1}\{(j,l)\}\right)$.
Let $\mathcal{F}\subseteq\left(\binom{[m]}{k}\right)$.  For a multiset $F\in \mathcal{F}$ and $i,j\in [m]$ and $s\in[k]$, define
$$S((i,s),j)(F) =\left\{
\begin{aligned}
F' & , & \text{ if  } (j,1)\notin F \text{ and } F'\notin \mathcal{F},\\
F &  , & \text{otherwise.}
\end{aligned}
\right.
$$   
Subsequently, define $S((i,s),j)(\mathcal{F}):=\{S((i,s),j)(F):F\in \mathcal{F}\}$.

Let $(\mathcal{F}, \mathcal{G})$ be a cross $t$-intersecting pair with $\mathcal{F},\mathcal{G}\subseteq\left( \binom{[m]}{k}\right)$.
A multiset $T$ is called a {\em $t$-kernel} of ($\mathcal{F}$, $\mathcal{G}$) if $|F\cap G\cap T|\geq t$ for any $F\in \mathcal{F}$ and $G\in \mathcal{G}$.
 Denote $\mathcal{K}(\mathcal{F},\mathcal{G})$ as the set of $t$-kernels of $(\mathcal{F}, \mathcal{G})$ such that $M(m,1)\subseteq T$ (or equivalently, $\{(1,1),(2,1), \dots, (m,1)\}\subseteq T$) for every $T\in\mathcal{K}(\mathcal{F},\mathcal{G})$. Clearly, $M(m,k)$ is a $t$-kernel, thus $\mathcal{K}(\mathcal{F},\mathcal{G})\not=\emptyset$. The following property holds.
\begin{lem}\label{PROP: tkernel}
 Let $(\mathcal{F}, \mathcal{G})$ be a cross $t$-intersecting pair in $\left( \binom{[m]}{k}\right)$. If $T$ is a $t$-kernel in $\mathcal{K}\left(\mathcal{F},\mathcal{G}\right)$, then $T \in \mathcal{K}\left(S\left((i,m(i,T)),j\right)(\mathcal{F}), S\left((i,m(i,T)),j\right)(\mathcal{G})\right)$ for any $i,j\in[m]$.
\end{lem}
\begin{proof}
According to the definition of down-compression, $S((i,m(i,T)),j)(\mathcal{F}$) and $S((i,m(i,T)),j)(\mathcal{G}$) are multiset families in $\left( \binom{[m]}{k}\right)$. 
Now we show that $S((i,m(i,T)),j)(\mathcal{F})$ and $S((i,m(i,T)),j)(\mathcal{G})$ are cross $t$-intersecting, and $T \in \mathcal{K}(S((i,m(i,T)),j)(\mathcal{F}),S((i,m(i,T)),j)(\mathcal{G}))$.
Choose $F\in \mathcal{F}$ and $G\in \mathcal{G}$ arbitrarily. We denote $S(F)=S((i,m(i,T)),j)(F)$ and $S(G)=S((i,m(i,T)),j)(G)$ and $s=m(i, T)$ for convenience in the following proof.

If $S(F)\in \mathcal{F}$ and $S(G)\in \mathcal{G}$, then $S(F)=F$ and $S(G)=G$ by the definition of down-compression. Therefore,
$$|S(F)\cap S(G)\cap T|=|F\cap G\cap T|\geq t,$$ we are done.

If $S(F)\notin \mathcal{F}$ and $S(G)\notin \mathcal{G}$, then according to  the definition of down-compression, 
we know that $(j,1)\notin F\cup G$,  however, $(j,1)\in S(F)\cap S(G)$. Since $s=m(i,T)$ with our choice, the multiplicity of element  $i$ decreases by  exactly one in $S(F)\cap S(G)\cap T$ compared to $F\cap G\cap T$. Note that $(j,1)\in T$.
Therefore,
$$|S(F)\cap S(G)\cap T|\geq \left|\left((F\cap G\cap T)\backslash \{(i,s)\}\right)\cup \{(j,1)\}\right|\geq t.$$

Finally, suppose that there exists exactly one of $F$ and $G$ that remains unchanged after the FGV-shifting. Without loss of generality, assume that  $S(F)\notin \mathcal{F}$ and $S(G)\in \mathcal{G}$, i.e., $S(F)=\left(F\setminus\cup_{t=s}^{m(i,F)}\{(i,t)\}\right) \bigcup \left(\cup_{l=1}^{m(i,F)-s+1}\{(j,l)\}\right)$ and $S(G)=G$. Since $S(G)=G$, we have $(j,1)\in G$, or $G'=\left(G\setminus\cup_{t=s}^{m(i,G)}\{(i,t)\}\right) \bigcup \left(\cup_{l=1}^{m(i,G)-s+1}\{(j,l)\}\right)\in\mathcal{G}$, or $m(i,T)> m(i,G)$. 
If $m(i,T)> m(i,G)$, then the multiplicity of element $i$ in $S(F)\cap S(G)$  remains unchanged compared to $F\cap G$. Therefore, 
$|S(F)\cap S(G)\cap T|\ge |F\cap G\cap T|\ge t$, we are done. Now assume $m(i,T)\le  m(i,G)$. Then the multiplicity of element  $i$ decreases by at most  one in $S(F)\cap S(G)\cap T$ compared to $F\cap G\cap T$.
If $(j,1)\in G$, note that $(j,1)\notin F$ but $(j,1)\in S(F)\cap G\cap T$, then $$|S(F)\cap S(G)\cap T|=|S(F)\cap G\cap T|\ge |F\cap G\cap T|-1+1\ge t.$$ 
Now assume $(j,1)\notin G$. Then $G'\in \mathcal{G}$ by the definition of down-compression, and the multiplicity of element $i$ in $S(F)\cap G$ is equal to its multiplicity in $F\cap G'$ (in fact, the multiplicity equals to $m(i,T)-1$). 
Obviously, $(j,1)\notin S(F)\cap G$ and $(j,1)\notin F\cap G'$. Therefore, $|S(F)\cap S(G)\cap T|=|S(F)\cap G\cap T|=|F\cap G'\cap T|\geq t$. 
\end{proof}


\begin{lem}\label{PROP: 4.3}
Let $(\mathcal{F}, \mathcal{G})$ be a cross $t$-intersecting pair in $\left(\binom{[m]}{k}\right)$ and let $T \in \mathcal{K}(\mathcal{F}, \mathcal{G})$ with $|T\backslash [m]|\geq 1$. For a fixed $i\in[m]$ with $m(i, T)\ge 2$, recursively define $S_{1,i}(\mathcal{F})=S((i,m(i,T)),1)(\mathcal{F})$, $S_{2,i}(\mathcal{F})=S((i,m(i,T)),2)(S_{1,i}(\mathcal{F}))$, $\dots$, and $S_{m,i}(\mathcal{F})=S((i,m(i,T)),m)(S_{m-1,i}(\mathcal{F}))$.
Similarly define $S_{m,i}(\mathcal{G})$. Then the following statements hold.

(1) $(S_{m,i}(\mathcal{F}), S_{m,i}(\mathcal{G}))$ is a cross $t$-intersecting pair with $|S_{m,i}(\mathcal{F})|=|\mathcal{F}|$ and $|S_{m,i}(\mathcal{G})|=|\mathcal{G}|$.
  
(2) $T\backslash\{(i, m(i,T))\}\in \mathcal{K}(S_{m,i}(\mathcal{F}), S_{m,i}(\mathcal{G}))$.

\end{lem}
\begin{proof}
(1) It is a direct corollary of Lemma~\ref{PROP: tkernel}.

For (2),  it is sufficient to prove that $|S_{m,i}(F)\cap S_{m,i}(G)\cap (T\backslash\{(i,m(i,T))\})|\geq t$ for any $F\in \mathcal{F}$ and $G\in \mathcal{G}$.
 According to Lemma~\ref{PROP: tkernel}, $T\in\mathcal{K}(S_{m,i}(\mathcal{F}), S_{m,i}(\mathcal{G}))$. If $(i,m(i,T))\notin S_{m,i}(F)\cap S_{m,i}(G)$, then $|S_{m,i}(F)\cap S_{m,i}(G)\cap (T\setminus \{(i, m(i,T))\})|=|S_{m,i}(F)\cap S_{m,i}(G)\cap T|\ge t$. 
 
 Now assume $(i,m(i,T))\in S_{m,i}(F)\cap S_{m,i}(G)$. Then $S_{m,i}(F)$ and $S_{m,i}(G)$ must remain unchanged after series of down-compression operations according to the definition of down-compression, i.e., $S_{m,i}(F)= F$ and $S_{m,i}(G)= G$. 
 Since $|F\cup G|\le 2k-t$ and $m(i,T)\geq2$, we have $|(F\cup G)\cap[m]|<2k-t$. Note that  $m\ge 2k-t$. By the pigeonhole principle,  there exists  $j\in [m]$ with $(j,1)\notin F\cup G$. However, $S_{m,i}(F)=F$, we must have  $S((i,m(i,T)),j)(F)\in \mathcal{F}$ according to the definition of down-compression. Thus $|S((i,m(i,T)),j)(F)\cap G\cap T|\geq t$. 
 Since $(i,m(i,T))\notin S((i,m(i,T)),j)(F)$, we have $(i,m(i,T))\notin S((i,m(i,T)),j)(F)\cap G\cap T$. 
Hence $|F\cap G\cap T|\ge |S((i,m(i,T)),j)(F)\cap G\cap T|+1\geq t+1$. Therefore, $|S_{m,i}(F)\cap S_{m,i}(G)\cap (T\backslash m(i,T))|=|F\cap G\cap (T\backslash m(i,T))|\geq t$.

\end{proof}

From Lemmas~\ref{PROP: tkernel} and \ref{PROP: 4.3}, we have the following corollary.
\begin{cor}\label{tkernel}
Let $(\mathcal{F}, \mathcal{G})$ be  a cross $t$-intersecting pair. Then there exists a bijection $f:\left(\binom{[m]}{k}\right) \rightarrow  \left(\binom{[m]}{k}\right)$ with the following property:

 (1)  $(f(\mathcal{F}), f(\mathcal{G}))$ is a cross $t$-intersecting pair;

 (2) $|\mathcal{F}|=|f(\mathcal{F})|$ and $|\mathcal{G}|=|f(\mathcal{G})|$;
 
 (3)  $M(m,1)$ is a $t$-kernel of  $(f(\mathcal{F}) ,f(\mathcal{G}))$.
 
\end{cor}
\begin{proof}
Choose a minimal $t$-kernel $T_0$ in $\mathcal{K}(\mathcal{F}, \mathcal{G})$. If $T_0=M(m,1)$, then take $f$ to be the identify mapping, and we are done. Now assume there exists $i_0$ such that $m(i_0, T_0)\ge 2$. 
Then apply the operation $S_{m,i_0}$ on $(\mathcal{F},\mathcal{G})$, we have an updated $t$-kernel $T_0\setminus\{(i,m(i_0,T_0))\}$ of $(S_{m,i_0}(\mathcal{F}), S_{m,i_0}(\mathcal{G}))$ according to Lemma \ref{PROP: 4.3}.
If $T_0\setminus\{(i,m(i_0,T_0))\}=M(m, 1)$, let $f=S_{m,i_0}$, then  we are done. If not, continue the operation $S_{m, i_1}$ on $(S_{m,i_0}(\mathcal{F}), S_{m,i_0}(\mathcal{G}))$ with $m(i_1, T_1)\ge 2$ for $T_1\in \mathcal{K}((S_{m,i_0}(\mathcal{F}), S_{m,i_0}(\mathcal{G})))$, then either we obtain a bijiection $f=S_{m, i_1}S_{m,i_0}$ and stop, or have an updated cross $t$-intersecting family with smaller $t$-kernel according to Lemma \ref{PROP: 4.3}. Since $|T_0|$ is finite, after applying finite times of the above operations, we finally have a bijection $f:\left(\binom{[m]}{k}\right) \rightarrow  \left(\binom{[m]}{k}\right)$ such that $M(m,1)\in \mathcal{K}(f(\mathcal{F}), f(\mathcal{G}))$.


\end{proof}

Now we can give the proof of Theorem~\ref{THM: t-intersecting}.

\begin{proof}

Let $(\mathcal{F}, \mathcal{G})$ be a maximum corss $t$-interscecting pair in $\left(\binom {[m]}{k}\right)$. 
By Corollary~\ref{tkernel}, we have a bijection $f: \left(\binom {[m]}{k}\right)\rightarrow \left(\binom {[m]}{k}\right)$ such that $M(m,1)$ $(=[m])$ is a $t$-kernel of  $(f(\mathcal{F}), f(\mathcal{G}))$.
Let $$\mathcal{X}_s=\{f(F)\cap[m]\,:\,f(F)\in f(\mathcal{F}) \text{ and } |f(F)\cap [m]|=s\}.$$ 
Then $t\le s\le k$.
Note that for a fixed set $X\in \mathcal{X}_s$,  the number of multisets $F\in \mathcal{F}$ with $F\cap [m]=X$ is at most $\binom{s+(k-s-1)}{k-s}=\binom{k-1}{k-s}$ (which equals to the number of the non-negative solutions of the  indeterminate equation $x_1+x_2+\dots+x_{k-s}=s$).
Therefore,  $|f(\mathcal{F})|\leq \sum_{s=t}^{k}|\mathcal{X}_s|\binom{k-1}{k-s}$.

Now consider a set family $\mathcal{F}'=\bigcup_{s=t}^k\mathcal{F}'_s$ on the ground set $[n]$, where $n=m+k-1$ and $\mathcal{F}'_s=\{F'\in\binom{[n]}{k}\, :\, F'\cap[m]\in\mathcal{X}_s\}$. 
Then $|\mathcal{F}'|=\sum_{s=t}^k|\mathcal{F}'_s|$. 
Note that for a fixed $X\in \mathcal{X}_s$,  the number of sets $F'\in \mathcal{F}'_s$ with $F'\cap [m]=X$ is  $\binom{k-1}{k-s}$. 
Therefore, $|\mathcal{F}'|\leq \sum_{s=t}^{k}|\mathcal{X}_s|\binom{k-1}{k-s}$,
and there exists a  injective mapping $\theta$ from $f(\mathcal{F})$ to $\mathcal{F}'$ such that $\theta(f(F))\cap[m]=f(F)\cap[m]$ (the mapping exists guaranteed by the cardinalities of $f(\mathcal{F})$ and $\mathcal{F'}$). 
Therefore, we have $|f(\mathcal{F})|\leq |\mathcal{F'}|$. 
Corresponding to $\mathcal{X}_s$ and $\mathcal{F}'$, we can define $\mathcal{Y}_s=\{f(G)\cap[m]\,:\,f(G)\in f(\mathcal{G}) \text{ with } |f(G)\cap [m]|=s\}$ and set family $\mathcal{G}'\in\binom{[n]}k$. With similar discussion, we have $|f(\mathcal{G})|\leq |\mathcal{G'}|$.

Now we claim that $\mathcal{F}'$ and $\mathcal{G}'\in \binom{[n]}{k}$ are  cross $t$-intersecting. Choose $F'\in\mathcal{F}'$ and $G'\in \mathcal{G}'$ arbitrarily. Then $|F'\cap G'|\geq |F'\cap G'\cap[m]|=|\theta^{-1}(F')\cap \theta^{-1}(G')\cap [m]|\ge t$ since $[m]=M(m,1)$ is a $t$-kernel of  $(f(\mathcal{F}), f(\mathcal{G}))$.

According to Corollary \ref{COR:JH},
\begin{eqnarray}\label{EQU: e1}
|\mathcal{F}|+|\mathcal{G}|=|f(\mathcal{F})|+|f(\mathcal{G})|\leq|\mathcal{F}'|+|\mathcal{G}'|\leq \binom{m+k-1}{k}-\sum\limits_{i=0}^{t-1}\binom{k}{i} \binom{m-1}{k-i}+1.     
\end{eqnarray}

In the following, we show that the extremal structure is unique up to isomorphism when equlity holds in (\ref{EQU: e1}). Note that $t\ge 2$. We first claim that the case $\{k,t\}=\{n-2,n-3\}$ does not happen. If it did, 
then we would have $k=n-2=m+k-3$, leading to $m=3$. Additionally, since $m\geq 2k-t$ and $k-t=1$, we must have $k=2$ and $t=1$. This contradicts to the initial condition $t\ge 2$. Therefore, when equality holds in (\ref{EQU: e1}), according to Corollary \ref{COR:JH}, we have unique extremal structure $\mathcal{F}'=\{[k]\}$ and $\mathcal{G}'=\{G'\in \binom{[m+k-1]}{k}  \,:\,G'\cap[k]\geq t\}$ up to isomorphism. 

Note that $$|\mathcal{F}|+|\mathcal{G}|=|f(\mathcal{F})|+|f(\mathcal{G})|=|\mathcal{F}'|+|\mathcal{G}'|.$$
Since  $|\mathcal{F}'|$=1, we have $|\mathcal{F}|=|f(\mathcal{F})|$=$|\mathcal{F}'|$=1. Let $\mathcal{F}=\{F\}$. Then $f(\mathcal{F})=\{f(F)\}$, implying that $\theta(f(F))=[k]$. Thus  $f(F)\cap [m]= [k]\cap[m]=[k]$, leading to $f(F)=[k]$. Therefore, $f(\mathcal{F})=\{[k]\}$. 
Since $(f(\mathcal{F}), f(\mathcal{G}))$ is a cross $t$-intersecting pair in $\left(\binom{[m]}{k}\right)$, we have $f(\mathcal{G})\subseteq\{G\in \left(\binom{[m]}{k}\right) \,:\,G\cap[k]\geq t\}$. 
Moreover, 
$$|f(\mathcal{G})|=|\mathcal{G}'|=\binom{m+k-1}{k}-\sum\limits_{i=0}^{t-1}\binom{k}{i} \binom{m-1}{k-i}=\left|\{G\in \left(\binom{[m]}{k}\right) \,:\,G\cap[k]\geq t\}\right|.$$
Therefore, we must have  $f(\mathcal{G})=\{G\in \left(\binom{[m]}{k}\right) \,:\,G\cap[k]\geq t\}$.

Denote $\mathcal{G}_k=\{G\in \left(\binom{[m]}{k}\right) \,:\,G\cap[k]\geq t\}$. It is sufficient to show that $\mathcal{F}=f(\mathcal{F})$ and $\mathcal{G}=f(\mathcal{G})$, i.e., $F=[k]$ and $\mathcal{G}=\mathcal{G}_k$ up to isomorphism. 
If $F\not=[k]$, then $|F\cap[m]|< k$. Since $(\mathcal{F}, \mathcal{G})$ is a maximum cross $t$-intersecting pair in $\left(\binom{[m]}{k}\right)$, we must have $\mathcal{G}=\{G\in \left(\binom{[m]}{k}\right) \,:\,|G\cap F|\geq t\}$. 
Note that $f(\mathcal{F}) (=\{[k]\})$ is obtained from $F$ through a series of down-compression operations.
Without loss of generality, we assume that $|F\cap[m]|= k-1$ (otherwise, we could repeatedly apply Lemma~\ref{PROP: 4.3} to $(\mathcal{F}, \mathcal{G})$ to increase the size of the support of $F$). 
Thus we may assume  $F=\{1*2,2*1,\dots,(k-1)*1\}$ up to isomorphism. Recall that $f(F)=[k]$ and $f(\mathcal{G})=\mathcal{G}_k$. 
Let $G_1=\{k-t+1,k-t+2,\dots, 2k-t\}$. Then $G_1\notin \mathcal{G}$ because $|G_1\cap F|=t-1<t$. Since each element of $G_1$ has multiplicity exactly one,  we have  $f(G_1)=G_1=\{k-t+1,k-t+2,\dots,2k-t\}$.  
Since $|G_1\cap [k]|=t$, we have $G_1\in \mathcal{G}_k=f(\mathcal{G})$, implying that $G_1=f^{-1}(f(G_1))\in\mathcal{G}$, a contradiction.

Therefore, $\mathcal{F}=\{[k]\}$ and $\mathcal{G}=\{G\in \left(\binom{[m]}{k}\right) \,:\,G\cap[k]\geq t\}$ are unique extremal structures  up to isomorphism.

\end{proof}

\section{Remarks and Discussions}	
In this paper, we characterize the extremal structures of cross $t$-intersecting multiset families for $t=1$ and $m\ge k+1$, and for $t\ge 2$ and $m\ge 2k-t$. Obviously, when $t=1$, we have a better lower bound of $m$ compared to $t\ge 2$. 
  However, for multiset intersecting families, $m$ can be less than $k$. Therefore, it will be interesting to consider the following question.
\begin{prob}
Let $\mathcal{F}$ and $\mathcal{G}$ are cross $t$-intersecting families in $\left(\binom{[m]}{k}\right)$. For $m<2k-t$ and $t>1$, what is the maximum of   $|\mathcal{F}|+|\mathcal{G}|?$
\end{prob}

 Denote $\binom{[n]_m}{k}=\{F\in \left(\binom{[n]}{k}\right)\,:\,m(i,F)\leq m, \text{ for all } i\in F\}$.
 A multiset family  $\mathcal{F} \subset \binom{[n]_m}{k}$  is called bounded if $m<k$, otherwise, $\mathcal{F}$ is considered unbounded. 
 In this article, we focus on the analysis of unbounded multisets.
In fact, Liao, Lv, Cao, and Lu~\cite{LL} presented an invariant of the Erd\H{o}s-Ko-Rado Theorem for the bounded multisets. 
\begin{thm}[\cite{LL}] 
    Let $m$, $k$, $n$ be integers. If $n\geq k+\lceil \frac{k}{m}\rceil$, then $|\mathcal{F}|\leq |\mathcal{S}_i|$, where $\mathcal{S}_i=\{F\in \binom{[n]_m}{k}\,:\,i\in F \}$ for some $i\in [n]$. When $n> k+\lceil \frac{k}{m}\rceil$, equality holds if and only if $\mathcal{F}=\mathcal{S}_i$
\end{thm}
Unfortunately, the techniques used in this paper are not valid for bounded multiset problems. 
We are interested in determining the maximum of the sum of cross $t$-intersecting bounded multiset families. 
\begin{prob}
Let $\mathcal{F}, \mathcal{G}\subset \binom{[n]_m}{k}$ be cross $t$-intersecting-families.
Determine the maximum of $|\mathcal{F}|+|\mathcal{G}|$.

\end{prob}

Another variant of cross $t$-intersecting problem involves maximizing the product of cross $t$-intersecting families. 
Let $\mathcal{F}$ and $\mathcal{G}$ be cross intersecting families.
Pyber \cite{PL} showed that 
$|\mathcal{F}||\mathcal{G}|\leq \binom{n-1}{k-1} \binom{n-1}{l-1}$ for $\mathcal{F}\subset \binom{[n]}{k}$ and $\mathcal{G}\subset \binom{[n]}{l}$ with $n\ge 2k+l-2$ and $k \geq l$.  Matsumoto and Tokushige \cite{MT} further demonstrated this upper bound holds for $n\geq 2k \geq 2l$. Tokushige in~\cite{PSMN} and Frankl,  Lee, Siggers and Tokushige in~\cite{TOKU} determined  the maximum value of the product of cross $t$-intersecting families $\mathcal{F}, \mathcal{G}\subset \binom{[n]}{k}$ for $n\ge 1,443(t+1)k$ (\cite{PSMN}) and $n\ge (t+1)k$ (\cite{TOKU}), respectively, conjecturing that the best lower bound for $n$ is $(t+1)(k-t+1)$. 
For cross $t$-intersecting multiset families, it is very interesting to consider the following question,
\begin{prob}
Let $\mathcal{F}$ and $\mathcal{G}$ be cross $t$-intersecting families in $\left(\binom{[n]}{k}\right)$.
Determine the maximum value of the product   $|\mathcal{F}||\mathcal{G}|$? Is it true that $|\mathcal{F}||\mathcal{G}|\le \binom{n+k-t-1}{k-t}^2?$
\end{prob}

\vspace{5pt}
\noindent{\bf Acknowledgements}:
This work was supported by the National Key Research and Development Program of China (2023YFA1010203), the National Natural Science Foundation of China (12471336), and the Innovation Program for Quantum Science and Technology (2021ZD0302902).

\hspace{-6.5mm}\textbf{Data Availability:}
Data sharing not applicable to this article as no datasets were generated or analysed during the current study.

\end{document}